\documentstyle[12pt]{article}

\newcommand{\F}{\noindent}

\newcommand{\MP}{\medskip}
\newcommand{\BP}{\bigskip}

\newcommand{\beq}{\begin{eqnarray}}
\newcommand{\ene}{\end{eqnarray}}
\newcommand{\beqn}{\begin{eqnarray*}}
\newcommand{\enen}{\end{eqnarray*}}

\newcommand{\aaa}{{\mbox{\bf a}}}

\newcommand{\qqqq}{{\mbox{\bf q}}}
\newcommand{\kk}{{\mbox{\bf k}}}
\newcommand{\xx}{{\mbox{\bf x}}}

\newcommand{\AAA}{{\bf A}}
\newcommand{\BBB}{{\bf B}}
\newcommand{\PP}{{\bf P}}

\newcommand{\eq}[1]{(\ref{#1})}

\Large

\begin{document}

\rightline{KIMS-2003-05-01}

\BP

\vskip12pt

\vskip8pt

\begin{center}
\Large
{\bf Is mathematics consistent?}
\vskip18pt

\normalsize
Hitoshi Kitada
\vskip2pt

Graduate School of Mathematical Sciences

University of Tokyo

Komaba, Meguro-ku, Tokyo 153-8914, Japan

e-mail: kitada@ms.u-tokyo.ac.jp
\vskip10pt
May 1, 2003
\end{center}

\BP

\leftskip24pt
\rightskip24pt

\small

\noindent
{\it Abstract}: A question is proposed whether or not set theory is consistent.

\leftskip0pt
\rightskip0pt

\vskip 24pt

\large

\F
We consider a formal set theory $S$, where we can develop a number theory. As no generality is lost, in the following we consider a number theory that can be regarded as a subsystem of $S$, and will call it $S^{(0)}$.
\BP

\F
{\it Definition} 1. 1) We assume that a G\"odel numbering of the system $S^{(0)}$ is given, and denote a formula with the G\"odel number $n$ by $A_n$.

\F
2) $\mbox{\AAA}^{(0)}(a,b)$ is a predicate meaning that ``$a$ is the G\"odel number of a formula $A$ with just one free variable (which we denote by $A(a)$), and $b$ is the G\"odel number of a proof of the formula $A(\aaa)$ in $S^{(0)}$," and $\mbox{\BBB}^{(0)}(a,c)$ is a predicate meaning that ``$a$ is the G\"odel number of a formula $A(a)$, and $c$ is the G\"odel number of a proof of the formula $\neg A(\aaa)$ in $S^{(0)}$." Here $\aaa$ denotes the formal natural number corresponding to an intuitive natural number $a$ of the meta level.
\MP

\F
{\it Definition} 2. Let $\mbox{\PP}(x_1,\cdots.x_n)$ be an intuitive-theoretic predicate. We say that $\mbox{\PP}(x_1,\cdots,x_n)$ is {\it numeralwise expressible} in the formal system $S^{(0)}$, if there is a formula $P(x_1,\cdots,x_n)$ with no free variables other than the distinct variables $x_1,\cdots,x_n$ such that, for each particular $n$-tuple of natural numbers $x_1,\cdots,x_n$, the following holds:
\MP

\F
i) if $\mbox{\PP}(x_1,\cdots,x_n)$ is true, then $\vdash P(\xx_1,\cdots,\xx_n).$
\MP

\F
and
\MP

\F
ii) if $\mbox{\PP}(x_1,\cdots,x_n)$ is false, then $\vdash \neg P(\xx_1,\cdots,\xx_n).$

\BP

\F
Here ``true" means ``provable on the meta level."

\BP

\F
{\bf Lemma 1}. There is a G\"odel numbering of the formal objects of the system $S^{(0)}$ such that the predicates $\mbox{\AAA}^{(0)}(a,b)$ and $\mbox{\BBB}^{(0)}(a,c)$ defined above are primitive recursive and hence numeralwise expressible in $S^{(0)}$ with the associated formulas $A^{(0)}(a,b)$ and $B^{(0)}(a,c)$. (See \cite{K}.)
\BP

\F
{\it Definition} 3. Let $q^{(0)}$ be the G\"odel number of a formula: 
$$
\forall b [\neg A^{(0)}(a,b)\vee \exists c(c\le b \hskip3pt\&\hskip2pt B^{(0)}(a,c))].
$$
Namely
\beq
A_{q^{(0)}}(a)=\forall b [\neg A^{(0)}(a,b)\vee \exists c(c\le b \hskip3pt\&\hskip2pt B^{(0)}(a,c))]\nonumber
\ene
In particular
\beq
A_{q^{(0)}}(\qqqq^{(0)})=\forall b [\neg A^{(0)}(\qqqq^{(0)},b)\vee \exists c(c\le b \hskip3pt\&\hskip2pt B^{(0)}(\qqqq^{(0)},c))]\nonumber
\ene

\BP

\ 

\BP

Assume that $S^{(0)}$ is consistent.
\BP

Suppose that
\beq
\vdash A_{q^{(0)}}(\qqqq^{(0)})\mbox{ in } S^{(0)},\nonumber
\ene
and let $k^{(0)}$ be the G\"odel number of the proof of $A_{q^{(0)}}(\qqqq^{(0)})$. Then by the numeralwise expressibility of $\AAA^{(0)}(a,b)$
\beq
\vdash A^{(0)}(\qqqq^{(0)},\kk^{(0)}).\label{Rosser1}
\ene
Under our hypothesis of consistency,
$$
\vdash A_{q^{(0)}}(\qqqq^{(0)})\mbox{ in } S^{(0)}
$$
implies
$$
\mbox{not } \vdash \neg A_{q^{(0)}}(\qqqq^{(0)})\mbox{ in } S^{(0)}.
$$
Hence, for any integer $\ell$, $\BBB^{(0)}(q^{(0)},\ell)$ is false. In particular, $\BBB^{(0)}(q^{(0)},0)$, $\cdots,$ $\BBB^{(0)}(q^{(0)},k^{(0)})$ are false. By virtue of the numeralwise expressibility of $\BBB^{(0)}(a,c)$, from these follows that
$$
\vdash \neg B^{(0)}(\qqqq^{(0)},0),\vdash \neg B^{(0)}(\qqqq^{(0)},1),\cdots, \vdash \neg B^{(0)}(\qqqq^{(0)},\kk^{(0)}).
$$
Hence
$$
\vdash \forall c (c\le \kk^{(0)} \supset \neg B^{(0)}(\qqqq^{(0)},c)).
$$
This together with $\vdash A^{(0)}(\qqqq^{(0)},\kk^{(0)})$ in \eq{Rosser1} gives
$$
\vdash \exists b [A^{(0)}(\qqqq^{(0)},b)\hskip3pt\&\hskip2pt\forall c(c\le b\supset\neg B^{(0)}(\qqqq^{(0)},c))].
$$
This is equivalent to
$$
\vdash \neg A_{q^{(0)}}(\qqqq^{(0)})\mbox{ in } S^{(0)}.
$$
A contradiction with our consistency hypothesis of $S^{(0)}$. Thus
$$
\mbox{not } \vdash A_{q^{(0)}}(\qqqq^{(0)})\mbox{ in } S^{(0)}.
$$

Reversely, suppose that
\beq
\vdash \neg A_{q^{(0)}}(\qqqq^{(0)})\mbox{ in } S^{(0)}.\nonumber
\ene
Then there is a G\"odel number $k^{(0)}$ of the proof of $\neg A_{q^{(0)}}(\qqqq^{(0)})$ in $S^{(0)}$, and we have
$$
\mbox{\BBB}^{(0)}(q^{(0)},k^{(0)}) \mbox{ is true.}
$$
Thus 
$$
\vdash B^{(0)}(\qqqq^{(0)},\kk^{(0)}),
$$
from which follows
\beq
\vdash \forall b[b\ge \kk^{(0)}\supset \exists c(c\le b\hskip3pt\&\hskip2pt B^{(0)}(\qqqq^{(0)},c))].\label{(1)}
\ene
As $\neg A_{q^{(0)}}(\qqqq^{(0)})$ is provable in $S^{(0)}$, there is no proof of $A_{q^{(0)}}(\qqqq^{(0)})$ in $S^{(0)}$ by our consistency assumption of $S^{(0)}$. Therefore
\beq
\vdash \neg A^{(0)}(\qqqq^{(0)},0), \vdash \neg A^{(1)}(\qqqq^{(1)},1),\cdots,\vdash \neg A^{(0)}(\qqqq^{(0)},\kk^{(0)}-1)\nonumber
\ene
hold. Thus
\beq
\vdash \forall b[b<\kk^{(0)}\supset \neg A^{(0)}(\qqqq^{(0)},b)]\label{(2)}.
\ene
Combining \eq{(1)} and \eq{(2)}, we obtain
$$
\vdash \forall b[\neg A^{(0)}(\qqqq^{(0)},b)\vee \exists c(c\le b\hskip3pt\&\hskip2pt B^{(0)}(\qqqq^{(0)},c))],
$$
which is
$$
\vdash A_{q^{(0)}}(\qqqq^{(0)}).
$$
A contradiction with our consistency assumption of $S^{(0)}$. Thus we have
\beq
\mbox{not }\vdash \neg A_{q^{(0)}}(\qqqq^{(0)})\mbox{ in } S^{(0)}.\nonumber
\ene

We have reproduced Rosser's form of G\"odel incompleteness theorem.
\BP

\F
{\bf Lemma 2}. Assume $S^{(0)}$ is consistent. Then neither $A_{q^{(0)}}(\qqqq^{(0)})$ nor $\neg A_{q^{(0)}}(\qqqq^{(0)})$ is provable in $S^{(0)}$.
\BP

Thus, we can add either one of $A_{q^{(0)}}(\qqqq^{(0)})$ or $\neg A_{q^{(0)}}(\qqqq^{(0)})$, which we will denote $A_{(0)}$ hereafter, as a new axiom of $S^{(0)}$ without introducing any contradiction. Namely, let $S^{(1)}$ be an extension of the formal system $S^{(0)}$ with an additional axiom $A_{(0)}$. Then by Lemma 2
\beq
S^{(1)} \mbox{ is consistent.}\label{consis}
\ene
We now extend definitions 1 and 3 to the extended system $S^{(1)}$ as follows with noting that the numeralwise expressibility of the predicates $\mbox{\AAA}^{(1)}(a,b)$ and $\mbox{\BBB}^{(1)}(a,c)$ defined below can be extended to the new system $S^{(1)}$ with {\it the same} G\"odel numbering {\it as} the one given in Lemma 1 for $S^{(0)}$.
\vskip12pt
\F
1) $\mbox{\AAA}^{(1)}(a,b)$ is a predicate meaning that ``$a$ is the G\"odel number of a formula $A(a)$, and $b$ is the G\"odel number of a proof of the formula $A(\aaa)$ in $S^{(1)}$," and $\mbox{\BBB}^{(1)}(a,c)$ is a predicate meaning that ``$a$ is the G\"odel number of a formula $A(a)$, and $c$ is the G\"odel number of a proof of the formula $\neg A(\aaa)$ in $S^{(1)}$."
\MP

\F
2) Let $q^{(1)}$ be the G\"odel number of a formula: 
$$
\forall b [\neg A^{(1)}(a,b)\vee \exists c(c\le b \hskip3pt\&\hskip2pt B^{(1)}(a,c))].
$$
\BP

\F
By the extended numeralwise expressibility, we have in the way similar to Lemma 2 by using the consistency \eq{consis} of $S^{(1)}$
\beq
\mbox{not } \vdash A_{q^{(1)}}(\qqqq^{(1)}) \quad\mbox {  and  } \quad\mbox{not } \vdash \neg A_{q^{(1)}}(\qqqq^{(1)}) \quad \mbox{ in } S^{(1)}.\label{pos1}\nonumber
\ene

Continuing the similar procedure, we get for any natural number $n(\ge 0)$ that
\beq
S^{(n)} \mbox{ is consistent,}\label{consis2}
\ene
\beq
\mbox{not } \vdash A_{q^{(n)}}(\qqqq^{(n)}) \quad\mbox {  and  }\quad \mbox{not } \vdash \neg A_{q^{(n)}}(\qqqq^{(n)}) \quad \mbox{ in } S^{(n)}.
\label{posn}
\ene

We now let $S^{(\omega)}$ the extended system of $S^{(0)}$ that includes all of the formulas $A_{(n)}(=A_{q^{(n)}}(\qqqq^{(n)})$ or $\neg A_{q^{(n)}}(\qqqq^{(n)}))$ as its axioms. By \eq{consis2} $S^{(\omega)}$ is consistent. We note that the formula $A_{(n)}$ is recursively defined if we have already constructed the system $S^{(n)}$. Moreover, if we let ${\tilde q}(n)$ be the G\"odel number of the formula $A_{(n)}$, we have ${\tilde q}(i)\ne{\tilde q}(j)$ for all natural numbers $i<j$ as $A_{(j)}$ is not provable in $S^{(i+1)}$ for $i<j$. Thus $\sup_{i\le n}{\tilde q}(i)$ goes to infinity as $n$ tends to infinity. Further we note that ${\tilde q}(n)$ is a recursive function of $n$. Then given a formula $A_r$ with G\"odel number $r$, restricting our attention to the formulas $A_{(n)}$ with $\tilde{q}(n)\le r$, we can determine in $S^{(\omega)}$ recursively if that given formula $A_r$ is an axiom of the form $A_{(n)}$ or not. Thus the addition of all $A_{(n)}$ retains the recursive definition of the following predicates $\mbox{\AAA}^{(\omega)}(a,b)$ and $\mbox{\BBB}^{(\omega)}(a,c)$ defined in the same way as above.
\BP

\F
$\mbox{\AAA}^{(\omega)}(a,b)$ is a predicate meaning that ``$a$ is the G\"odel number of a formula $A(a)$, and $b$ is the G\"odel number of a proof of the formula $A(\aaa)$ in $S^{(\omega)}$," and $\mbox{\BBB}^{(\omega)}(a,c)$ is a predicate meaning that ``$a$ is the G\"odel number of a formula $A(a)$, and $c$ is the G\"odel number of a proof of the formula $\neg A(\aaa)$ in $S^{(\omega)}$."
\BP

\F
Then we see that the predicates $\mbox{\AAA}^{(\omega)}(a,b)$ and $\mbox{\BBB}^{(\omega)}(a,c)$ are numeralwise expressible in $S^{(\omega)}$ and the G\"odel number $q^{(\omega)}$ of the formula: 
$$
\forall b [\neg A^{(\omega)}(a,b)\vee \exists c(c\le b \hskip3pt\&\hskip2pt B^{(\omega)}(a,c))],
$$
denoted by $A_{q^{(\omega)}}(a)$, is well-defined.

As before, we continue the similar procedure, transfinite inductively. In this process, from the nature of our extension procedure, at each step $\alpha$ where we construct the $\alpha$-th consistent system $S^{(\alpha)}$ from the preceding systems $S^{(\gamma)}$ with $\gamma<\alpha$, the predicates $\AAA^{(\alpha)}(a,b)$ and $\BBB^{(\alpha)}(a,c)$ must be recursively defined from the preceding predicates $\AAA^{(\gamma)}(a,b)$ and $\BBB^{(\gamma)}(a,c)$ $(\gamma<\alpha)$ so as for $S^{(\alpha)}$ to be further extended with retaining consistency. For this to hold, it is necessary and sufficient that the ordinal $\alpha$ is a recursive ordinal. 

Is there any ordinal that is not recursive? A function $F(x)$ is called recursive if it has the form:
$$
F(x)=G(x,F|x),
$$
where $F|x$ is a restriction of $F$ to a domain $x$, and $G$ is a given function. Consider the formula for an ordinal $x$:
$$
x=\{y | y\in x\}.
$$
This meets the above requirement of the recursive definition of the ordinal $x$, although this is tautological and may not be considered a definition usually. But if we see the structure that it contains the domain $x$ and by using that domain only it defines $x$ itself, it is not so unreasonable to think that there is no nonrecursive ordinal.

There is, however, a possibility (\cite{M}) that the condition whether or not a nonrecursive ordinal exists in ZFC is independent of the axioms of ZFC. In that case we have two alternatives.
\MP

\F
Case i) There is no nonrecursive ordinal, and hence all ordinals are recursive.
\MP

\F
In this case, the extension of the system $S^{(\alpha)}$ above is always possible. Thus we can extend $S^{(\alpha)}$ indefinitely forever. However, in this process, we cannot reach the step where the number of added axioms is the cardinality $\aleph_1$ of the first uncountable ordinal, as the number of added axioms is at most countable by the nature of formal system. Thus there must be a least countable ordinal $\beta$ such that the already constructed consistent system $S^{(\beta)}$ is not extendable with retaining consistency. This contradicts the unlimited extendibility stated above, and we have a contradiction. Insofar as we assume that every ordinal is recursive, the only possibility remaining is to conclude that set theory is inconsistent.
\MP

\F
Case ii) There is a nonrecursive ordinal, thus there is a least nonrecursive ordinal $\omega_1$ usually called Church-Kleene ordinal (\cite{F}, \cite{T}).
\MP

\F
In this case the above extension of $S^{(\alpha)}$ is possible if and only if $\alpha<\omega_1$. We note that $\omega_1$ is a limit ordinal. For if it is a successor of an ordinal $\delta$, then $\delta<\omega_1$ is recursive, hence so is $\omega_1=\delta+1$, a contradiction with the nonrecursiveness of $\omega_1$. Therefore we can construct, in the same way as that for $S^{(\omega)}$, a consistent system $S^{(\omega_1)}$, which cannot be extended further with retaining consistency by the nonrecursiveness of $\omega_1$.

On the other hand, as we have seen in the discussion of case i), there must be a least countable ordinal $\beta$ such that the already constructed consistent system $S^{(\beta)}$ is not extendable with retaining consistency. Since $\beta$ is the least ordinal such that $S^{(\beta)}$ is not extendable, for any $\alpha<\beta$ the system $S^{(\alpha)}$ is consistently extendable. Whence by the reasoning above about the recursiveness of $\alpha$ with which $S^{(\alpha)}$ is consistently extendable, $\alpha$ is recursive and we have $\alpha<\omega_1$ if $\alpha<\beta$. Thus
\beq
\beta\le \omega_1.\label{beta}
\ene
Reversely, when $\alpha<\omega_1$, $\alpha$ is a recursive ordinal. Thus by the same reasoning as above about the recursiveness of $\alpha$, $S^{(\alpha)}$ is consistently extendable. Therefore $\alpha<\beta$ if  $\alpha<\omega_1$. This and \eq{beta} give
$$
\beta=\omega_1.
$$

\BP

Summarizing, we have proved

\BP

\F
{\bf Theorem}. Assume that $S^{(0)}$ is consistent. Suppose that the condition whether or not there is a nonrecursive ordinal is independent of the axioms of ZFC. Then there are the following two alternatives:
\MP

\F
i) There is no nonrecursive ordinal, and hence all ordinals are recursive.
\begin{quotation}
\F
In this case set theory is inconsistent.
\end{quotation}

\F
ii) There is a nonrecursive ordinal, thus there is a least countable nonrecursive ordinal $\beta$.
\begin{quotation}
\F
In this case the corresponding system $S^{(\beta)}$ is consistent and cannot be extended further with retaining consistency.
\end{quotation}
\BP

We remark that this is a metamathematical theorem. 

Thus the inconsistency in i) of this theorem does not give any proof in ZFC of the existence of nonrecursive ordinal. To know whether a nonrecursive ordinal exists or not, we need a proof in ZFC or if such a statement is independent of the axioms of ZFC, we need to add an axiom that determines which the case is. In the latter case, the above theorem shows a direction in which the extended ZFC can be consistent if the original ZFC is consistent.

Further, as the above theorem is a metamathematical theorem, even if there is no nonrecursive ordinal, the case i) of the theorem does not yield that set theory is inconsistent in the sense that we can find a concrete inconsistent proposition like Russell's paradox inside the set theory. Rather it would be said that we may not find such an inconsistent proposition insofar as we work inside the set theory ZFC. Thus this theorem should not be interpreted as stating any concrete inconsistency of set theory.

\BP

\end{document}